\documentclass[%
 jmp,
 amsmath,amssymb,
notitlepage]{revtex4-1}

\usepackage{graphicx}
\usepackage{dcolumn}
\usepackage{bm}

\usepackage[utf8]{inputenc}
\usepackage[T1]{fontenc}
\usepackage{mathptmx}
\newcommand{\beq}{\begin{equation}}
\newcommand{\eeq}{\end{equation}}
\newcommand{\bea}{\begin{eqnarray}}
\newcommand{\eea}{\end{eqnarray}}
\newcommand{\bwd}{\begin{widetext}}
\newcommand{\ewd}{\end{widetext}}

\begin{document}

\preprint{AIP/123-QED}

\title[Sample title]{A Proof of Riemann Hypothesis}

\author{Tao Liu}
\affiliation{
School of Science, Southwest University of Science and Technology, Mianyang, Sichuan 621010, China\\
State Key Laboratory of Environment-friendly Energy Materials, Southwest University of Science and Technology, 59 Qinglong Road, Mianyang, Sichuan 621010, China
}%

\author{Juhao Wu}
\affiliation{%
Stanford University, Stanford, California 94309, USA
}%

\date{\today}

\begin{abstract}
The meromorphic function $W(s)$ introduced in the Riemann-Zeta function $\zeta(s) = W(s) \zeta(1-s)$ maps the line of $s = 1/2 + it$ onto the unit circle in $W$-space. $|W(s)| = 0$ gives the trivial zeroes of the Riemann-Zeta function $\zeta(s)$. In the range: $0 < |W(s)| \neq 1$, $\zeta(s)$ does not have nontrivial zeroes. $|W(s)|=1$ is the necessary condition for the nontrivial zeros of the Riemann-Zeta function. Writing $s = \sigma + it$, in the range: $0 \leq \sigma \leq 1$, but $\sigma \neq 1/2$, even if $|W(s)|=1$, the Riemann-Zeta function $\zeta(s)$ is non-zero. Based on these arguments, the nontrivial zeros of the Riemann-Zeta function $\zeta(s)$ can only be on the $s = 1/2 + it$ critical line. Therefore a proof of the Riemann Hypothesis is presented.
\end{abstract}

\maketitle

\section{Riemann Hypothesis}
Let us briefly revisit the Riemann Hypothesis. Recall that the Riemann-Zeta function can be defined via Riemann's functional equation:
    \beq\label{zeta}
    \zeta(s) = 2^s\pi^{s-1}\sin\left(\frac{\pi s}2\right)\Gamma(1-s)\zeta(1-s),
    \eeq
where $\Gamma(1-s)$ is the analytical continuation of the factorial and $s = \sigma + i t$, with $\sigma \in R$ and $t \in R$ both being real number. Notice that $s = - 2 n$ ($n = 1, 2, \cdots, \infty$) are the trivial zeroes of $\zeta(s)$ and $s = 1$ is its pole. The Riemann Hypothesis states on the necessary condition of the nontrivial zeroes of the Riemann-Zeta function as: ``all the nontrivial zeroes of $\zeta(s)$ is on the line: $s = 1/2 + i t$'' \cite{Rieman59}. Up to now, it is proven that the nontrivial zeroes of $\zeta(s)$ can only be in the range: $0\leq \sigma \leq 1$ \cite{Hadamard96, Vallee-Poussin96}

\section{The Meromorphic Function $W(s)$}\label{Riemann}

\subsection{Introduction of the meromorphic function $W(s)$}
Based on Riemann's functional equation (\ref{zeta}), we can introduce a meromorphic function:
    \beq
    W(s) = 2^s \pi^{s-1} \sin\left(\frac{\pi s}2\right)\Gamma(1-s),
    \eeq
so that Eq. (\ref{zeta}) is rewritten as:
    \beq\label{zetawz}
    \zeta(s) = W(s) \zeta(1-s).
    \eeq
The distribution of the nontrivial zeros of the Riemann-Zeta function is closely related to the properties of the meromorphic function $W(s)$. In the following, let us discuss the properties of $W(s)$.

\subsection{The properties of the meromorphic function $W(s)$}

\subsubsection{The reflection symmetry of $W(s)$}\label{property1} 
\begin{figure}[!htp]
 \centering
 \includegraphics[width=6cm,angle=0]{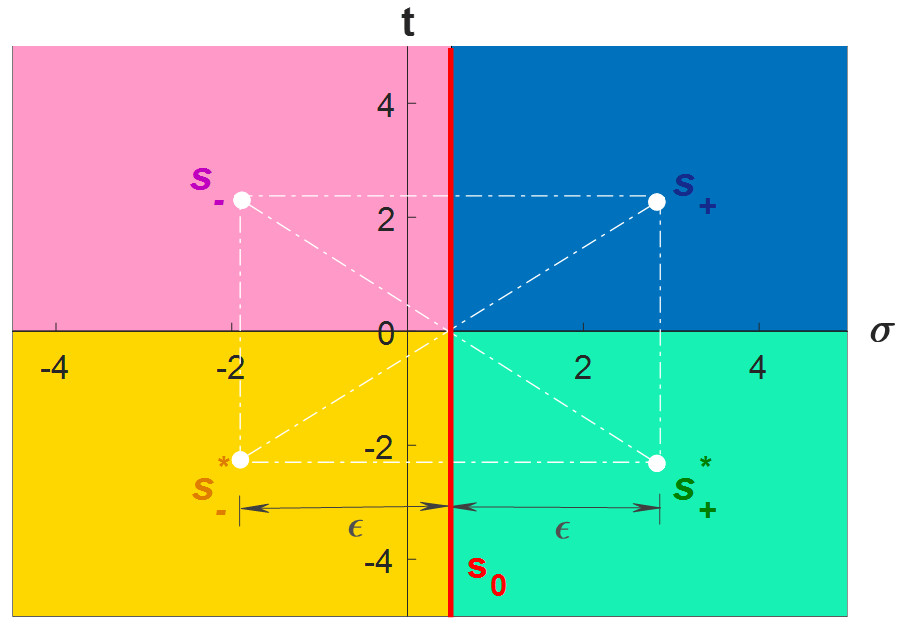}
 \caption{\label{symmetricPair}
 In the $s$-complex plane, $s_{+}$ and $s_{-}$ are mirror symmetric with respect to $s_0$; $s^{\ast}_{+}$ and $s^{\ast}_{-}$are mirror symmetric with respect to $s_0$.
 }
\end{figure}

In the $s$-complex space, for arbitrary $\varepsilon \in R$, setting $s_0 = 1/2 + it$, then the pair: $s_{\pm} = s_0 \pm \varepsilon$ is a mirror symmetric pair with respect to $s_0$ as shown in Fig. \ref{symmetricPair}. The complex conjugates of $s_+$ and $s_-$ are noted as $s^{\ast}_+$ and $s^{\ast}_-$. Under the $W(s)$ map, $s_{\pm}$ and their mirror reflected complex conjugate $s^{\ast}_{\mp}$ are reciprocal pair.

    \beq\label{ww1}
    W\left(s_{\pm}\right)W\left(s^{\ast}_{\mp}\right) = 1.
    \eeq

\noindent {\it Proof}

Because:
$$\sin\left(\frac{\pi s}2\right) = \frac{\pi}{\Gamma(1-s/2)\Gamma(s/2)},$$ and
$$\frac{\Gamma(1-s)}{\Gamma(1-s/2)} = \frac{2^{-s} \Gamma\left(\frac{1-s}2\right)}{\sqrt{\pi}},$$
we can rewrite $W(s)$ as:
    \beq\label{wgamma}
    W(s)=\pi^{s-\frac12}\frac{\Gamma\left(\frac{1-s}2\right)}{\Gamma\left(\frac s2\right)}.
    \eeq
On the complex plane of $s$, for arbitrary real number $\varepsilon \in R$,
because:

    $$
    W(s_{+})
    = e^{\varepsilon + i t} \frac{\Gamma\left(\frac14 - \frac{\varepsilon}2 - \frac{it}2\right)}{\Gamma\left(\frac 14 + \frac{\varepsilon}2 + \frac{it}2\right)},
    $$
and
    $$
    W(s^{\ast}_{-})
    =
    e^{-\varepsilon - i t}\frac{\Gamma\left(\frac14 + \frac{\varepsilon}2 + \frac{it}2\right)}{\Gamma\left(\frac 14 - \frac{\varepsilon}2 - \frac{it}2\right)} = \frac 1{W(s_+)},
    $$
which is to say that in the $W$-space, $W(s^{\ast}_-)$ must be equal to $1/W(s_{+})$, which is the reciprocal of $W(s_+)$; we have:
    $$
    W(s_+)W(s^{\ast}_-) = W(s_+) \frac1{W(s_+)} = 1.
    $$
Similarly, we can prove that in the $W$-space, $W(s^{\ast}_+)$ must be equal to $1/W(s_{-})$, which is the reciprocal of $W(s_-)$, {\it i.e.}, $W(s^{\ast}_+)W(s_-) = 1$. Therefore, property \ref{property1} is proven.

\subsubsection{$W(s)$ has trivial zeros at $s = \{- 2n | n = 0, 1,2, \cdots, \infty\}$ and pole $s = \{2 n + 1 | n = 0, 1, 2, \cdots, \infty \}$. $W(s)$ at these trivial zeros and at these poles are reciprocal pairs: $W(-2n)W(1 + 2n) = 1$.}\label{reflective}\hfill

\noindent {\it Proof}

First of all, it is obvious that:
$$
W(-2n) = 2^{-2n}\pi^{-2n-1}\sin(-n\pi)\Gamma(1+2n) = 0
$$
for ($n=0,1,2,\cdots$),
and
$$
W(2n+1) = 2^{2n+1}\pi^{2n}\sin\left(n\pi + \frac{\pi}2 \right)\Gamma(-2n) = \infty$$
for ($n=0,1,2,\cdots$).

Next, setting $\varepsilon = \frac12 + 2n$, and $t = 0$, we have $s_+ = 2n + 1$, and $s_- = -2n$. According to reflection symmetry relation as in Eq. (\ref{ww1}), we have:
    \beq
    W(2n+1)W(-2n) = 1.
    \eeq
This then serves as the proof for property \ref{reflective}.

\subsubsection{$W(s)$ is not zero at $s = 2n$ ($n=1,2,\cdots$).}\label{wnz}\hfill

\noindent {\it Proof}

Let $\varepsilon = - \frac12 + 2 n$, and $t = 0$, then $s_+ = 2n$, and $s_- = 1 - 2n$. According to the reflection symmetry as in Eq. (\ref{ww1}), we have:

$$
W(2n)W(1-2n) = 1,
$$
{\it i.e.},
    \beq
    W(2n)
    =
    \frac1{W(1-2n)}
    =
    \frac1{2^{1-2n}\pi^{-2n}\sin\left(\frac{\pi}2 - n\pi\right)\Gamma(2n)}
    =
    \frac{(2\pi)^{2n}}{2(-1)^n(2n-1)!}\neq 0.
    \nonumber
    \eeq
This is a proof for property \ref{wnz}.

\subsubsection{Monotonicity of the absolute value $|W(s)|$ of $W(s)$. \newline \indent Excluding the zeros and the poles, the absolute value $|W(s)|$ of the map $W(s)$ is a monotonic function of $t$ except when $\sigma = 1/2$: \newline \indent 1. in the range $0 < t < +\infty$, when $\sigma > 1/2$, $|W(s)|$ monotonically decreases with the increase of $t$; when $\sigma < 1/2$, $|W(s)|$ monotonically increases with the increase of $t$. \newline \indent 2. in the range $-\infty < t < 0$, when $\sigma > 1/2$, $|W(s)|$ monotonically increases with the increase of $t$; when $\sigma < 1/2$, $|W(s)|$ monotonically decreases with the increase of $t$.}\label{mono}\hfill

\noindent {\it Proof}

The fact that $W(s^{\ast}) = 2^{s^{\ast}}\pi^{s^{\ast}-1}\sin\left(\frac{\pi s^{\ast}}2\right)\Gamma\left(1-s^{\ast}\right) = \overline W(s)$, excluding the zeros and the poles of $W(s)$, because $|W(s)|^2 = W(s)W(s^{\ast})$, we only need to take the first-order derivative with respect to $t$ on both sides of the equation:
    \beq\label{dwdt}
    2|W(s)|  \frac d{dt} |W(s)| = \frac d{dt} \left[W(s)W\left(s^{\ast} \right)\right].
    \eeq
The right hand side of Eq. (\ref{dwdt}) gives:
    \beq\label{dwwdt}
    \frac d{dt} \left[W(s)W\left(s^{\ast} \right)\right]
    =
    \frac{i|W(s)|^2}2 \left[\Psi\left(\frac{\sigma -it}2\right) - \Psi\left(\frac{\sigma + it} 2\right)
    -
    \Psi\left(\frac{1-\sigma-it}2\right) + \Psi\left(\frac{1-\sigma + it}2\right)\right]
    \eeq
where $\Psi(z)$ is the digamma function. Therefore, except for the zeros and the poles of $W(s)$, according to Eq. (\ref{dwdt}) and Eq. (\ref{dwwdt}), we have:
    \beq\label{dabswdt}
    \frac d{dt} |W(s)|
    =
    \frac{i|W(s)|}4 \left[\Psi\left(\frac{\sigma -it}2\right) - \Psi\left(\frac{\sigma + it} 2\right) \right.
    -
    \left.\Psi\left(\frac{1-\sigma-it}2\right) + \Psi\left(\frac{1-\sigma + it}2\right)\right]
    \eeq
where the series expression of $\Psi(z)$ is:
    \beq\label{Psigamma}
    \Psi(\sigma+it)=-\gamma + \sum_{n=1}^{\infty}\frac{\sigma + it - 1}{n(n + \sigma + it - 1)}
    \eeq
with $\gamma$ being the Euler constant. Further noticing that:
    $$\Psi\left(\frac{\sigma -it}2\right) - \Psi\left(\frac{\sigma + it}2\right) = \sum_{n=1}^{\infty}\frac{-4it}{|it+2n+\sigma-2|^2},$$
and
    $$-\Psi\left(\frac{1-\sigma-it}2\right) + \Psi\left(\frac{1-\sigma +it}2\right) = \sum_{n=1}^{\infty}\frac{4it}{|it+2n-\sigma-1|^2},$$
so for a more explicit expression showing its being positive or negative, Eq. (\ref{dabswdt}) can be rewritten as:
    \beq\label{dwdtpos}
    \frac d{dt} |W(s)|=
    \left(\frac12-\sigma\right)t|W(s)| \sum_{n=1}^{\infty}\frac{8\left(n-\frac34\right)}{|it+2n+\sigma-2|^2 |it+2n-\sigma-1|^2}.
    \eeq
With the summation being positively defined, and $|W(s)|>0$, therefore:

1) For arbitrary non-zero $t$, when and only when $\sigma =1/2$, we have $\frac d{dt} |W(s)|=0$. This is to say that $\sigma =1/2$ is the only zero of $\frac d{dt} |W(s)|$.

2) When $\sigma \neq 1/2$, the value $\frac d{dt} |W(s)|$ being positive or negative is uniquely determined by the coefficient $(1/2-\sigma)t$ in Eq. (\ref{dwdtpos}).

\noindent In the range: $0<t<\infty$,

  when $\sigma <1/2$, $\frac d{dt} |W(s)|>0$, therefore $|W(s)|$ monotonically increases with the increase of $t$;

  when $\sigma >1/2$, $\frac d{dt} |W(s)|<0$, therefore $|W(s)|$ monotonically decreases with the increase of $t$.

\noindent In the range: $-\infty<t<0$,

  when $\sigma <1/2$, $\frac d{dt} |W(s)|<0$, therefore $|W(s)|$ monotonically decreases with the increase of $t$;

  when $\sigma >1/2$, $\frac d{dt} |W(s)|>0$, therefore $|W(s)|$ monotonically increases with the increase of $t$.

\noindent So, property \ref{mono} is proven.

\subsubsection{In the range: $0 \leq \sigma \leq 1$, but $\sigma  \neq 1/2$, the $t$ satisfying $|W(s)|=1$ is bounded: $2\pi < |t| < \kappa$.}\label{bounded}\hfill

\noindent {\it Proof}

1) Because $|W(s)|=|W(s^{\ast})|$, the $t$ satisfying $|W(s)|=1$ is always symmetric with respect to $\pm t$. So, let us only study $t \geq 0$ case.

2) Because $W(s_+) W(s_-^{\ast}) = 1$, we have $|W(s_+)||W(s_-)| = 1$, {\it i.e.}, when $s_{\pm} = \frac 12 \pm \varepsilon + it$, $|W(s_+)|$ and $|W(s_-)|$ have reflection symmetry.

3) In the range: $0 < \sigma < 1$,

when $t=0$, $|W(\sigma)| = W(\sigma)$, it is easy to prove $\frac {dW(\sigma)}{d\sigma} > 0$ (please refer to Appendix \ref{appendpos}), so $W(\sigma)$ monotonically increases with the increase of $\sigma$. Furthermore:
$$0 \leq W(\sigma) \leq1, \hspace{0.25cm}   (0\leq \sigma \leq 1/2);$$
and
$$1 \leq W(\sigma) \leq \infty, \hspace{0.25cm} (1/2 \leq \sigma \leq 1).$$

When $t$ increases from $0$, based on the monotonicity of $|W(s)|$ with respect to $t$ and the reflection symmetry of $\left|W\left(s_{+}\right)\right|\left|W\left(s_{-}\right)\right| = 1$, we know:

for a given $\sigma$ being ($0\leq \sigma <1/2$), $|W(\sigma + it)|$ must monotonically increase to being larger than $1$ from being less than $1$; likewise, for a given $\sigma$ being ($1/2<\sigma \leq 1$), $|W(\sigma +it)|$ must monotonically decrease to being less than $1$ from being larger than $1$. Therefore, there always exists $t_1 < t_2$, so that $\sigma_{\pm} = \frac12 \pm \varepsilon$ with $0< \varepsilon \leq 1/2$:
$$|W(\sigma_{-} + it_1 )| < 1 < |W(\sigma_{-} + it_2 )|,\hspace{0.25cm} (0 \leq \sigma_{-} < 1/2 ).$$ Due to reflection symmetry, $|W(\sigma_{-} + i t_1)||W(\sigma_{+} + i t_1)| = 1$, $|W(\sigma_{-} + i t_2)||W(\sigma_{+} + i t_2)| = 1$, one must have:
$$|W(\sigma_{+} + it_2 )| < 1 < |W(\sigma_{+} + it_1 )|,\hspace{0.25cm} (1/2 < \sigma_{+} \leq 1).$$

4) $t_2 = 2.01 \pi$

When $t_2 = 2.01 \pi$, $|W(0 + i t_2)| = 1.0025 > 1$, $\left|W\left(\frac12 +  i t_2 \right)\right| = 1$, $|W(1 + i t_2)| = \frac1{1.0025} < 1$. It is easy to prove $\left. \frac {d|W(\sigma + it)|}{d\sigma}\right|_{t = t_2} < 0$ (please refer to Appendix \ref{appendneg}), therefore $|W(\sigma + i t_2)|$ monotonically decreases with the increase of $\sigma$. Therefore, we have:
$$|W(\sigma + it_2 )| > 1, \hspace{0.25cm} (0 \leq \sigma < 1/2 );$$
and
$$|W(\sigma + it_2 )| < 1, \hspace{0.25cm} (1/2 < \sigma \leq 1).$$

This is to say that when $t$ changes from $0$ to $t_2$, in the range $0 \leq \sigma < 1/2$, the value $|W(\sigma + i t)|$, starting from being smaller than $1$, always monotonically increases to be larger than $1$; while in the range $1/2 < \sigma \leq 1$, the value $|W(\sigma + i t)|$, starting from being larger than $1$, always monotonically decreases to be smaller than $1$. Therefore, when $\sigma \neq 1/2$, the $t$ satisfying $|W(\sigma + it)| = 1$ has an up limit, {i.e.}, $t$ won't be larger than $\kappa \equiv t_2 = 2.01 \pi$. In the range $0 \leq \sigma \leq 1$, for arbitrary $\sigma \neq 1/2$, due to the fact that the $t$ satisfying $|W(s)| = 1$ is always reflection symmetric with respect to $\pm t$, $t$ has a low limit $-\kappa$, and an up limit $\kappa$.

5) $t_1 = 2 \pi$

When $t_1 = 2 \pi$, $|W(0 + i t_1)| = 0.9999999991 < 1$, $\left| W\left( \frac12 + i t_1 \right)\right| = 1$, $|W(1 + i t_1)| = \frac1{0.9999999991} > 1$. It is easy to prove (please refer to Appendix \ref{largsmal}) that:
$$|W(\sigma + it_1)| < 1\hspace{0.25cm} (0\leq \sigma < 1/2),$$
and
$$|W(\sigma + it_1)| > 1\hspace{0.25cm} (1/2 < \sigma \leq 1).$$
So in the range $0 \leq \sigma \leq 1$, and $\sigma \neq 1/2$, the $t$ satisfying $|W(\sigma + it)| = 1$ can only be in the range $t_1 < |t| < t_2$. That is to say, the $t$ satisfying $|W(s)| = 1$ is bounded:
    \beq
    2\pi < |t| < \kappa.
    \eeq

\begin{figure}[!htp]
 \centering
 \includegraphics[width=6cm,angle=0]{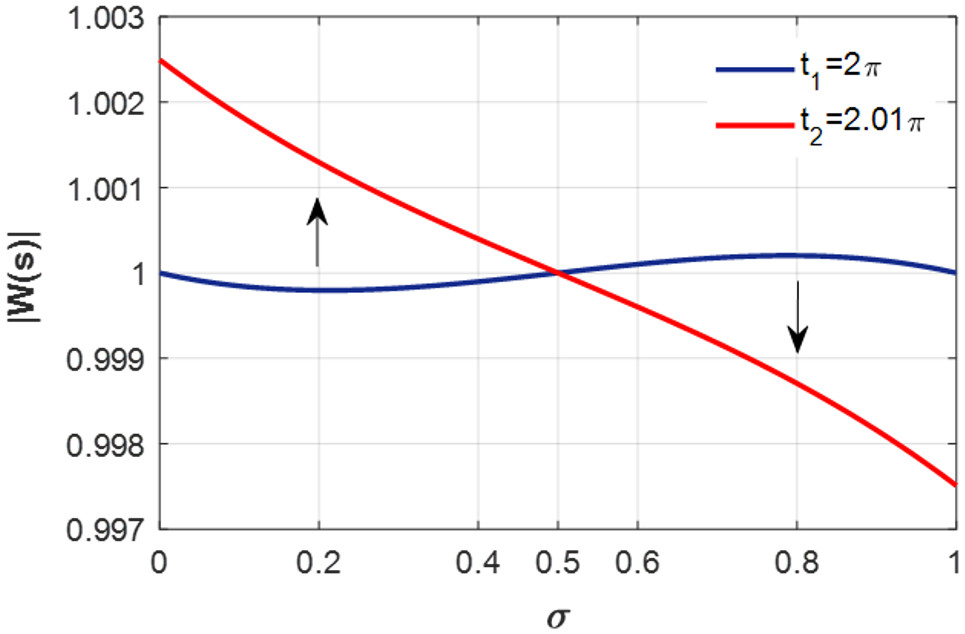}
 \caption{\label{wssig}
 When $t = 2\pi$, and $2.01\pi$, the continuous evolution of $|W(\sigma +it)|$ in the range of $0 \leq \sigma \leq 1$. In the plot: $\uparrow$ stands for the monotonic increase direction of $|W(s)|$ with $t$; while $\downarrow$ stands for the monotonic decrease direction of $|W(s)|$ with $t$.
 }
\end{figure}

Just for illustration, with $t_1 = 2 \pi$ and $t_2 = 2.01 \pi$, the continuous evolution of $|W(\sigma + it)|$ with $\sigma$ is shown in Fig. \ref{wssig}. 

\noindent So, property \ref{bounded} is proven.

\section{The distribution of the nontrivial zeros of the Riemann-Zeta function}\label{Proof}
Let us take the absolute value of the Riemann Equation in (\ref{zetawz}) to have:
    \beq
    |\zeta(s)| = |W(s)| |\zeta(1-s)|.
    \eeq
In the following, we will first prove a few lemmas of the Riemann-Zeta function $\zeta(s)$ based on the properties of the meromorphic function $W(s)$.

{\bf Lemma 1: the set of the zeroes of $W(s)=0$: $s=\{-2n | n=0,1,2,\cdots,\infty\}$ contains all the trivial zeroes of Riemann-Zeta function $\zeta(s)$. There is no nontrivial zero in this set.}

\noindent {\it Proof}

Notice that the set of all the trivial zeroes of $\zeta(s)$ is: $s=\{ -2n | n=1,2,\cdots,\infty\}$, which contains all the elements of the set of the zeroes of $W(s)$ except $n=0$. However, when $n=0$, we have $\zeta(0)=-1/2 \neq 0$. Therefore, Lemma $1$ is proven.

{\bf Lemma 2: Riemann-Zeta function $\zeta(s)$ and $\zeta(1-s)$ do not have nontrivial zeroes in the range: $0<|W(s)|\neq1$.}

\noindent {\it Proof}

{\bf 1) If $|W(s)| \neq 1$, then it must be true that $|\zeta(s)| \neq |\zeta(1 - s)|$.}

Assuming that $|\zeta(s)| = |\zeta(1-s)|$, then the necessary condition for the equation $|\zeta(s)| = |\zeta(1-s)|$ being equivalent to the Riemann equation $|\zeta(s)| = |W(s)| |\zeta(1-s)|$ is $|W(s)| = 1$ (please refer to Appendix \ref{equiv}). Then this is in conflict of $|W(s)| \neq 1$. Therefore, the assumption could not hold and it has to be true that $|\zeta(s)| \neq |\zeta(1-s)|$.

{\bf 2) When $|W(s)| > 0$, for the Riemann-Zeta function satisfying $|\zeta(s)| \neq |\zeta(1-s)|$, it is impossible to have $|\zeta(s)| = 0$.}

Assuming that $|\zeta(s)| = 0$, then due to $|\zeta(s)| = |W(s)| |\zeta(1-s)|$, it must be true that $|\zeta(1-s)| = 0$, so that $|\zeta(s)| = |\zeta(1-s)|$. Then this is in conflict with the statement that $|\zeta(s)| \neq |\zeta(1-s)|$. Therefore, the assumption $\zeta(s) = 0$ can not hold.

{\bf 3) When $|W(s)| > 0$, for the Riemann-Zeta function satisfying $|\zeta(s)| \neq |\zeta(1-s)|$, $|\zeta(1 - s)| = 0$ can only be valid on the trivial zeroes.}

Assuming that $|\zeta(1 - s)| = 0$, then due to $|\zeta(s)| = |W(s)| |\zeta(1-s)|$, one would have $|\zeta(s)| = 0$ with the only exception of $|W(s)| = \infty$, {\it i.e.}, $s$ is a pole of $|W(s)|$. However, $|\zeta(s)| = 0$ is in conflict with $|\zeta(s)| \neq |\zeta(1-s)|$. So, for $|\zeta(1 - s)| = 0$ and also $|\zeta(s)| \neq |\zeta(1-s)|$, the only possibility is $|W(s)| = \infty$, {\it i.e.}, $s$ is the pole of $|W(s)|$. Now, the poles of $|W(s)|$ are $s = 2 n + 1$, where $|\zeta(1-s)| = |\zeta(-2n)| = 0$, {\it i.e.}, these are just the trivial zeroes.

Based on the above 1), 2), and 3), we know that the Riemann-Zeta function does not have nontrivial zeroes for $0 < |W(s)| \neq 1$. So Lemma $2$ is proven.

{\bf Lemma 3: $|W(s)|=1$ is the necessary condition for the nontrivial zeroes of the Riemann-Zeta function $\zeta(s)$.}\label{w1necess}

\noindent {\it Proof}

Based on Lemma $1$: the set of zeroes from $W(s)=0$, {\it i.e.}, $s=\{-2n | n=0,1,2,\cdots,\infty\}$ does not contain nontrivial zeroes of the Riemann-Zeta function $\zeta(s)$.

Based on Lemma $2$: In the range: $0<|W(s)|\neq1$, there is no nontrivial zeroes of the Riemann-Zeta function $\zeta(s)$.

Therefore, the nontrivial zeroes of the Riemann-Zeta function can only be on $|W(s)|=1$. Indeed, when $|W(s)|=1$, we must have:
    \beq
    |\zeta(s)|=|\zeta(1-s)|.
    \eeq
If $|\zeta(s)|=0$, then $|\zeta(1-s)|=0$. However, $|\zeta(s)|=|\zeta(1-s)|$ does not guarantee $\zeta(s)$ to be zero. Therefore, $|W(s)|=1$ is only the necessary condition of the nontrivial zeroes of the Riemann-Zeta function $\zeta(s)$.

So Lemma $3$ is proven.

{\bf Corollary 1: for $s=1/2+it$ being the nontrivial zeroes of $\zeta(s)$, the necessary condition for $W(s)$ is that $|W(s)|=1$.}

\begin{figure*}[!htp]
 \centering
 \includegraphics[width=9cm,angle=0]{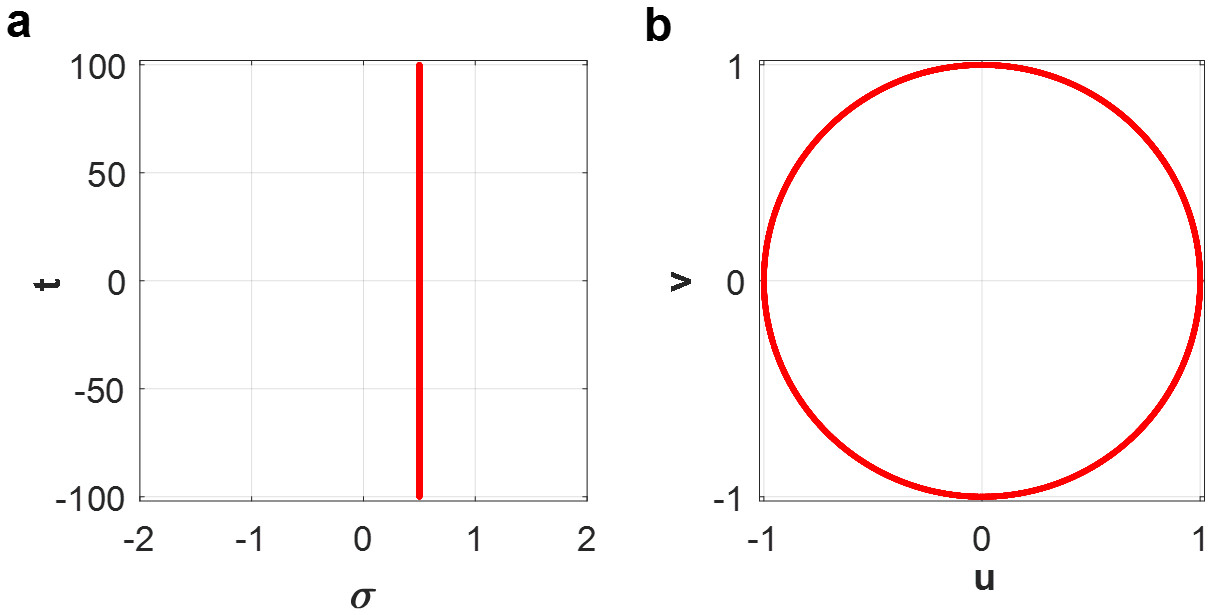}
 \caption{\label{ws}
 Writing $W(s)=u(\sigma,t)+iv(\sigma,t)$, the meromorphic function $W(s)$ maps the straight line $s=1/2+it$ in $s$-space (a) onto the unit circle in $W$-space (b). In the Figure, there are $8,000$ points in the range: $t=-100\cdots100$.
 }
\end{figure*}

Inserting $s = \frac12 + i t$ and $s^{\ast} = \frac12 - i t$ into Eq. (\ref{wgamma}), we have:
    \beq
    W(s){\overline W(s)} = \pi^{s+s^{\ast}-1} \left.\frac{\Gamma\left(\frac{1-s}2\right)\Gamma\left(\frac{1-s^{\ast}}2\right)}{\Gamma\left(\frac s2\right) \Gamma\left(\frac {s^{\ast}}2\right)}\right|_{s=\frac12+it,s^{\ast}=\frac12-it} = 1.
    \eeq
So we have: $|W(s)|=1$. The geometric illustration is that $W(s)$ maps the straight line $s=1/2+it$ in $s$-space to the unit circle in the $W$-space as in Fig. \ref{ws} (a) and (b).

{\bf Lemma 4: Riemann-Zeta function does not have nontrivial zeroes in the range $0 \leq \sigma < 1/2$ and $1/2 < \sigma \leq 1$; and $|t| > \kappa$ and $|t|<2\pi$.}

\noindent {\it Proof}

According to Lemma $3$, $|W(s)| = 1$ is the necessary condition for the Riemann-Zeta function to have nontrivial zeroes. On the other hand, according to the property \ref{bounded} of the meromorphic function: for $|W(s)| = 1$ but $\sigma \neq 1/2$, $t$ is bounded, {\it i.e.}, $2\pi < |t| < \kappa$ with $\kappa = 2.01\pi$. Therefore, in the range $|t| > \kappa$ and $|t| < 2\pi$, we have $|W(s)| \neq 1$. Hence, except $\sigma = 1/2$, Riemann-Zeta function does not have nontrivial zeros in the range $0 \leq \sigma \leq 1$ and $|t| < 2\pi$ and $|t| > \kappa$.
So Lemma $4$ is proven.

{\bf Lemma 5: Riemann-Zeta function does not have nontrivial zeroes in the range $0 \leq \sigma < 1/2$ and $1/2 < \sigma \leq 1$; and $2\pi < |t| < \kappa$.}

\noindent {\it Proof}

1) In the range $0 \leq \sigma < 1/2$, in order to prove that the Riemann-Zeta function $\zeta(s) \neq 0$ and $\zeta(1-s) \neq 0$ when $|W(s)|$ monotonically increases from $0 < |W(s)| < 1$ to $|W(s)| = 1$, we only need to prove that when $|W(s)| = 1$, it is not in the form of $\frac 00$.

According to Eq. (\ref{zetawz}) and Eq. (\ref{wgamma}), we have:
    \beq\label{wzetazeta}
    |W(s)| = \frac{|\zeta(s)|}{|\zeta(1-s)|} = \pi^{\sigma - \frac12}\frac{\left| \Gamma\left(\frac{1-s}2\right)\right|}{\left| \Gamma\left(\frac s2 \right)\right|}.
    \eeq
Notice that when $t > 0$ (please refer to Appendix \ref{dgammadt}):
    \beq
    \frac d{dt}\left| \Gamma\left( \frac {1 - s} 2\right)\right| = -t \sum^{\infty}_{n=1}\frac{\left| \Gamma\left( \frac {1 - s} 2\right)\right|}{|it+ 2n - \sigma -1|^2} < 0,
    \eeq
    \beq
    \frac d{dt}\left| \Gamma\left( \frac s 2\right)\right| = -t \sum^{\infty}_{n=1}\frac{\left| \Gamma\left( \frac s 2\right)\right|}{|it+ 2n + \sigma -2|^2} < 0.
    \eeq
So, both $\left| \Gamma\left( \frac {1 - s} 2\right)\right|$ and $\left| \Gamma\left( \frac s 2\right)\right|$ are continuous functions which monotonically decrease with increasing $t$. For arbitrary $0 \leq \sigma < \frac12$, $\pi^{\sigma - 1/2}$ is bounded and not equal to zero. Therefore, when and only when $t \rightarrow \infty$, we have $\left| \Gamma\left( \frac {1 - s} 2\right)\right| \rightarrow 0$, and $\left| \Gamma\left( \frac s 2\right)\right| \rightarrow 0$. According to Eq. (\ref{wzetazeta}), we have:
    \beq\label{lim00}
    \lim_{t\rightarrow \infty} \frac{|\zeta(s)|}{|\zeta(1-s)|} = \lim_{t\rightarrow \infty} \pi^{\sigma - \frac12}\frac{\left| \Gamma\left(\frac{1-s}2\right)\right|}{\left| \Gamma\left(\frac s2 \right)\right|} \rightarrow \frac00.
    \eeq
Because when $t > 0 $, both $\left| \Gamma\left(\frac{1-s}2\right)\right|$ and $\left| \Gamma\left(\frac s2\right)\right|$ are continuous function with no singular point, when $\sigma \neq \frac12$, the ratio of the absolute values of the Riemann-Zeta functions $\zeta(s)$ and $\zeta(1-s)$: $\frac{|\zeta(s)|}{|\zeta(1-s)|}$ continuously monotonically increase with increasing $t$, when and only when $t \rightarrow \infty$, can it be in the form of $\frac00$. Therefore, logically for a certain arbitrary finite $t > 0$, $\frac{|\zeta(s)|}{|\zeta(1-s)|}$ can not be in the form of $\frac00$ (please refer to Appendix \ref{appendzetazetano00}). In particular, when $t$ continuously changes from $2\pi$ to $2.01 \pi$, the $t$ satisfying $|W(s)| = 1$, according to the property \ref{bounded}, has an up limit, $t < \kappa$. In this case, according to Eq. (\ref{wzetazeta}) and Eq. (\ref{lim00}), we have:
    \beq
    1 = |W(s)| = \frac{|\zeta(s)|}{|\zeta(1-s)|} = \pi^{\sigma - \frac12}\frac{\left| \Gamma\left(\frac{1-s}2\right)\right|}{\left| \Gamma\left(\frac s2 \right)\right|} \nrightarrow \frac00.
    \eeq
Therefore in the range:  $0 \leq \sigma < 1/2$ and $2\pi < |t| < \kappa$, the Riemann-Zeta function satisfying $|W(s)| = 1$ does not have nontrivial zero.

2) Due to the fact that $|W(s)|$ is reflection symmetric with respect to $\sigma = 1/2$, in the range:  $1/2 < \sigma \leq 1 $, the Riemann-Zeta function satisfying $|W(s)| = 1$ does not have nontrivial zero, either.

\begin{figure*}[!htp]
 \centering
 \includegraphics[width=6cm,angle=0]{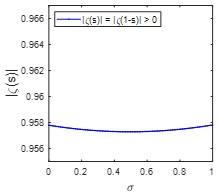}
 \caption{\label{lemma5plot}
 A numerical example of the absolute values of the Riemann-Zeta functions $|\zeta(s)|$ and $|\zeta(1-s)|$ in the range $0 \leq \sigma \leq 1$, but $\sigma \neq \frac 12$; while satisfying $|W(s)| = 1$. Notice that we do not intentionally exclude the point of $\sigma = \frac 12$ on the curve, since here we just want to give an illustration.
 }
\end{figure*}

3) In the range $0 \leq \sigma \leq 1$, but $\sigma \neq \frac12$, numerical calculations as shown in Fig. \ref{lemma5plot} of the absolute value of the Riemann-Zeta functions satisfying $|W(s)| = 1$ agree with the theoretical predictions in 1) and 2).

So Lemma $5$ is proven.

Now based on Lemma $1$, Lemma $2$, Lemma $3$, Lemma $4$, Lemma $5$, and Corollary $1$, we finally reach the conclusion that the nontrivial zeroes of the Riemann-Zeta function can only locate on the $s = 1/2 + it$ critical line. The Riemann Hypothesis is then proven.

\section{Discussion}\label{Discussion}
In this paper, we introduce a meromorphic function $W(s)$ via the Riemann functional equation: $\zeta(s) \equiv  W(s) \zeta(1-s)$. By studying the properties of the absolute value of this meromorphic function $|W(s)|$, we find that the necessary condition for Riemann-Zeta function to have nontrivial zero is $|W(s)| = 1$. With this necessary condition, we not only prove the Riemann Hypothesis in the range $0 \leq \sigma \leq 1$, but also we give a direct explanation why Riemann-Zeta function does not have nontrivial zero in the range $\sigma \leq 0$ and $\sigma \geq 1$ on the complex plane based on the zeros and singular points of $|W(s)|$.

1) The zeros of $|W(s)|$ correspond to the singular points of $\left| \Gamma\left( \frac s2 \right)\right| = |\Gamma(-n)|$. While $\left| \Gamma\left( \frac {1-s}2 \right)\right| = |\Gamma\left(\frac{2n+1}2\right)|$ are not singular.

2) Both the non-singular $\left| \Gamma\left( \frac s2 \right)\right|$ and the non-singular $\left| \Gamma\left( \frac {1-s}2 \right)\right|$ monotonically decrease, when $t$ monotonically increases near the zeros of $|W(s)|$. When and only when $t \rightarrow \infty$, we have $\left| \Gamma\left( \frac {1-s}2 \right)\right| \rightarrow 0$, and $\left| \Gamma\left( \frac s2 \right)\right| \rightarrow 0$, {\it i.e.},
$$
    \lim_{t\rightarrow \infty} \frac{|\zeta(s)|}{|\zeta(1-s)|} = \lim_{t\rightarrow \infty} \pi^{\sigma - \frac12}\frac{\left| \Gamma\left(\frac{1-s}2\right)\right|}{\left| \Gamma\left(\frac s2 \right)\right|} \rightarrow \frac00.
$$

3) Near its zeros: ($-2n$), ($0 < |W(s)| < 1$), during the monotonic increase of $t$, $|W(s)|$ will complete all the $\{s = -2n + \varepsilon + it\}$ points satisfying $|W(s)| = 1$ in the $\sigma < \frac 12$ complex plane. In this case:
$$
    1 = |W(s)| = \frac{|\zeta(s)|}{|\zeta(1-s)|} = \pi^{-2n+\varepsilon - \frac12} \frac{\left|\Gamma\left( n + \frac{1-\varepsilon}2 - \frac{it}2\right)\right|}{\left|\Gamma\left( -n + \frac{\varepsilon}2 + \frac{it}2\right)\right|} \nrightarrow \frac00,
$$
and $t$ is bounded.

This is to say that the $\frac{|\zeta(s)|}{|\zeta(1-s)|} = 1$ satisfying $|W(s)| = 1$ is not in the form of $\frac00$. Therefore, $|\zeta(s)|$ and $|\zeta(1-s)|$ both do not have nontrivial zeros in the $\sigma < \frac12$ complex plane. Furthermore, based on the reflection symmetry of $W(s)$ with respect to $\sigma = \frac12$, we readily know that the Riemann-Zeta function does not have nontrivial zero in the $\sigma > \frac12$ complex plane.

The above statement agrees with the known conclusion that the nontrivial zeroes of $\zeta(s)$ can only be in the range: $0\leq \sigma \leq 1$ \cite{Hadamard96, Vallee-Poussin96}; and the nontrivial zeroes of the Riemann-Zeta function can only locate on the $s = 1/2 + it$ critical line \cite{Rieman59}.

In summary, a proof of the Riemann Hypothesis is presented in this paper.

\appendix

\section{Proof of $\left.\frac{d|W(s)|}{d\sigma}\right|_{t=0} > 0$}\label{appendpos}
Because:

    \beq\label{dabswdsigma}
    \frac{d|W(s)|}{d\sigma}
    =
    \frac{|W(\sigma)|}4 \left[4 \ln(\pi) - \Psi\left(\frac s2 \right) - \Psi\left( \frac{s^{\ast}}2 \right)
    -
    \Psi\left(\frac{1 - s}2 \right) - \Psi\left(\frac{ 1- s^{\ast}}2\right)\right]
    \eeq
and making use of:
$$
\Psi(s)=-\gamma + \sum_{n=1}^{\infty}\frac{s - 1}{n(n + s - 1)},
$$
when $t = 0$, we have:
    \beq\label{apds}
    \Psi\left(\frac s2 \right) + \Psi\left( \frac{s^{\ast}}2 \right) = -2\gamma - \sum_{n=1}^{\infty}\frac{2(2 - \sigma)}{n(2n + \sigma - 2)},
    \eeq
and
    \beq\label{apd1ms}
    \Psi\left(\frac {1-s}2 \right) + \Psi\left( \frac{1- s^{\ast}}2 \right) = -2\gamma - \sum_{n=1}^{\infty}\frac{2(1 + \sigma)}{n(2n - \sigma - 1)}.
    \eeq
Now plugging Eq. (\ref{apds}) and Eq. (\ref{apd1ms}) into Eq. (\ref{dabswdsigma}), we have:
    \beq\label{dabswdsigt}
    \left.\frac{d|W(s)|}{d\sigma}\right|_{t=0}
    =
    \frac{|W(\sigma)|}4 \left\{4 \ln(\pi) + 4 \gamma + \sum_{n=1}^{\infty}\left[\frac{2(2 - \sigma)}{n(2n + \sigma - 2)}
    +
    \frac{2(1 + \sigma)}{n(2n - \sigma - 1)}\right]\right\}.
    \eeq
Now, it is clear that each term in the infinite summation in Eq. (\ref{dabswdsigt}) is positive when $0 < \sigma < 1$, therefore it is true that $\left. \frac{d|W(s)|}{d\sigma} \right|_{t=0} > 0$.

\section{Proof of $\left.\frac{d|W(s)|}{d\sigma}\right|_{t=2.01\pi} < 0$}\label{appendneg}
Setting:

    \beq\label{gsigt}
    G(\sigma, t)
    =
    4 \ln(\pi) - \Psi\left(\frac s2 \right) - \Psi\left( \frac{s^{\ast}}2 \right) - \Psi\left(\frac{1 - s}2 \right)
    -
    \Psi\left(\frac{ 1- s^{\ast}}2\right)
    \eeq
then Eq. (\ref{dabswdsigma}) can be rewritten as:
    \beq
    \frac{d|W(s)|}{d\sigma} = \frac{|W(s)|}4 G(\sigma, t).
    \eeq
In order to prove $\left.\frac{d|W(s)|}{d\sigma}\right|_{t=2.01\pi} < 0$, one needs only to prove the maximum value of $G(\sigma, 2.01\pi)$ satisfies: $\max[G(\sigma, 2.01\pi)] < 0$.

{\bf 1) when $\sigma = 1/2$ and $t = 2.01 \pi$, $G(1/2, 2.01 \pi)$ is the maximum value}

The first order and second order derivatives of $G(\sigma, t)$ with respect to $\sigma$ are:
    \beq\label{fisrtder}
    \frac{dG(\sigma, t)}{d\sigma}
    =
    \frac12 \left[ -\Psi\left(1,\frac s2\right) - \Psi\left( 1, \frac{s^{\ast}}2 \right) + \Psi\left( 1, \frac{1-s}2 \right)
    +
    \Psi\left( 1, \frac{1 - s^{\ast}}2 \right) \right],
    \eeq
and
    \beq
    \frac{d^2G(\sigma, t)}{d\sigma^2}
    =
    \frac14 \left[ -\Psi\left(2,\frac s2\right) - \Psi\left( 2, \frac{s^{\ast}}2 \right) - \Psi\left( 2, \frac{1-s}2 \right)
    -
    \Psi\left( 2, \frac{1 - s^{\ast}}2 \right) \right],
    \eeq
where $\Psi(1,z) = \frac{d\Psi(z)}{dz}$, $\Psi(2,z) = \frac{d^2\Psi(z)}{dz^2}$. Because when $\sigma = \frac12$, $1 -s = s^{\ast}$:
$$
\Psi\left( 1, \frac{1 - s}2\right) = \Psi\left( 1, \frac{s^{\ast}}2 \right),
$$
and
$$
\Psi\left( 1, \frac{1 - s^{\ast}}2\right) = \Psi\left( 1, \frac s2\right).
$$
Hence, we have:
$$
\left. \frac{dG(\sigma, t)}{d\sigma} \right|_{\sigma = \frac 12} = 0,
$$
while
$$
\left. \frac{d^2G(\sigma, t)}{d\sigma^2} \right|_{\sigma = \frac 12, t = 2.01 \pi} = -1.009542407 < 0.
$$
Therefore, $G\left(\frac 12, 2.01 \pi \right)$ is the maximum value, and
    \beq
    G\left( \frac12, 2.01 \pi \right) = - 0.015751728 < 0.
    \eeq

{\bf 2) when $t = 2.01 \pi$, $G\left( \frac 12, 2.01 \pi \right)$ is the maximum value in the range $0 \leq \sigma \leq 1$}

Setting
    \beq
    F(\sigma, t) = \left[ \Psi\left(\frac s2\right) + \Psi\left(\frac{s^{\ast}}2\right) + \Psi\left(\frac {1-s}2\right) + \Psi\left(\frac{1 - s^{\ast}}2\right)\right],
    \eeq
then
    \beq\label{gf}
    G(\sigma, t) = 4 \ln(\pi) - F(\sigma, t).
    \eeq
Because in the range $0 \leq \sigma \leq 1$, and $t \neq 0$, $F(\sigma, t)$ is an analytical function, there exists any order of derivative of $F(\sigma, t)$ with respect to $\sigma$:
    \beq\label{dnfdsign}
    \frac{d^n}{d\sigma^n}F(\sigma, t) = \frac1{2^n}\left[ \Psi\left(n, \frac s2\right) + (-1)^n \Psi\left(n, \frac {1-s^{\ast}}2\right) + \Psi\left(n, \frac {s^{\ast}}2\right) + (-1)^n \Psi\left(n, \frac {1-s}2\right) \right].
    \eeq
When $\sigma = \frac12$, we have: $s = 1 - s^{\ast} = \frac12 + it$, and $s^{\ast} = 1 - s = \frac12 - it$, plugging into Eq. (\ref{dnfdsign}), we have:
    \beq
    \left. \frac{d^n}{d\sigma^n}F(\sigma, t) \right|_{\sigma = 1/2} = \left\{
    \begin{array}{cc}
    0,& (n=2k+1) \\
    \frac1{2^{2k -1}}\left[ \Psi\left( 2k, \frac14 + \frac{it}2\right) + \Psi\left( 2k, \frac14 - \frac{it}2\right) \right], & (n = 2k)
    \end{array}
    \right..
    \eeq
Let us Taylor expand $F(\sigma, t)$ at $\sigma = \frac12$, we have:
    \beq\label{fsigt}
    F(\sigma, t) = \left. \sum^{\infty}_{n = 0} \frac1{n!} \frac{d^n}{d\sigma^n}F(\sigma, t) \right|_{\sigma = \frac12} \left( \sigma - \frac12 \right)^n = \sum^{\infty}_{k=0} \frac{\Psi\left( 2k, \frac14 + \frac{it}2 \right) + \Psi\left( 2k, \frac14 - \frac{it}2 \right)}{(2k)!2^{2k-1}} \left( \sigma - \frac12 \right)^{2k}.
    \eeq
Therefore,
    \beq\label{dfdsig}
    \frac d{d\sigma}F(\sigma, t) = \sum^{\infty}_{k=1} \frac{2k\left[\Psi\left( 2k, \frac14 + \frac{it}2 \right) + \Psi\left( 2k, \frac14 - \frac{it}2 \right)\right]}{(2k)!2^{2k-1}} \left( \sigma - \frac12 \right)^{2k-1} = \sum^{\infty}_{k=1} \frac{\Psi\left( 2k, \frac14 + \frac{it}2 \right) + \Psi\left( 2k, \frac14 - \frac{it}2 \right)}{(2k-1)!2^{2k-1}} \left( \sigma - \frac12 \right)^{2k-1}.
    \eeq
Setting $\sigma - \frac12 = \delta$, according to Eq. (\ref{gf}) and Eq. (\ref{fsigt}), we get the Taylor expansion of $G(\sigma, t)$:
    \beq\label{gsigtdel}
    G(\sigma, t) = 4\ln(\pi) - \sum^{\infty}_{k=0} \frac{\Psi\left(2k, \frac14 + \frac{it}2\right) + \Psi\left(2k, \frac14 - \frac{it}2 \right)}{2^{2k-1}(2k)!}\delta^{2k}.
    \eeq
Notice that:
$$
\frac d{d\sigma}G(\sigma, t) = - \frac d{d\sigma} F(\sigma, t),
$$
and plugging into Eq. (\ref{dfdsig}), we readily get the Taylor expansion of $\frac d{d\sigma}G(\sigma, t)$:
    \beq
    \frac d{d\sigma}G(\sigma, t) = -\sum^{\infty}_{k=1}\frac{\Psi\left(2k, \frac14 + \frac{it}2\right) + \Psi\left(2k, \frac14 - \frac{it}2\right)}{2^{2k-1}(2k-1)!}\delta^{2k-1}.
    \eeq

In particular, setting $t = 2.01 \pi$, we have:
    \beq
    \frac{dG(\sigma, 2.01\pi)}{d\sigma} = c_1\delta + c_3\delta^3 + c_5\delta^5 +O(\delta^7),
    \eeq
where $c_{2k-1} = -\frac{\Psi\left(2k, \frac14 + \frac{2.01\pi t}2\right) + \Psi\left(2k, \frac14 - \frac{2.01\pi t}2\right)}{(2k-1)!2^{2k-1}}$; in particular,
$$c_1 = -1.0095424 \times 10^{-1},$$
$$c_3 = +2.5705715 \times 10^{-3},$$
and
$$c_5 = -6.6264219 \times 10^{-5}.$$
For a sufficiently large $N$, $|c_5| > \sum^N_{k=4}|c_{2k-1}| \rightarrow 1.824122120 \times 10^{-6}$, i.e., neglecting the terms with order higher than $O(\delta^7)$ does not change the conclusion on $\frac{dG(\sigma, 2.01\pi)}{d\sigma}$ being positive or negative.

When $\delta < 0$: $c_1 \delta > 0$, $c_3 \delta^3 < 0$, $c_5 \delta^5 > 0$, and $c_1 \delta + c_3 \delta^3 > 0$, so at the left side of the zero $\left(-\frac12 \leq \delta < 0\right)$, $\left. \frac d{d\sigma} G(\sigma, t) \right|_{t=2.01\pi} > 0$, so $G(\sigma, 2.01\pi)$ monotonically increases with $\sigma$.

When $\delta > 0$: $c_1 \delta < 0$, $c_3 \delta^3 > 0$, $c_5 \delta^5 < 0$, and $c_1 \delta + c_3 \delta^3 < 0$, so at the right side of the zero $\left(0 < \delta \leq \frac12\right)$, $\left. \frac d{d\sigma} G(\sigma, t) \right|_{t=2.01\pi} < 0$, so $G(\sigma, 2.01\pi)$ monotonically decreases with $\sigma$.

So, $G(1/2, 2.01\pi)$ must be the maximum in the range $0 \leq \sigma \leq 1$, {\it i.e.}, $\max[G(\sigma, 2.01\pi)] = G(1/2, 2.01\pi)$.

\section{Proof of $|W(\sigma_{-} + 2 \pi i)| < 1$; $|W(\sigma_{+} + 2 \pi i)| > 1$ with $\left(\sigma_{\pm} = \frac12 \pm \varepsilon\right)$}\label{largsmal}

{\it Proof}

Based on the reflection symmetry, we have $|W(\sigma_{-} + 2 \pi i)||W(\sigma_{+} + 2 \pi i)| = 1$. To prove $|W(\sigma_{+} + 2 \pi i)| > 1$, we only need to prove $|W(\sigma_{-} + 2\pi i)| < 1$.

{\bf 1) $G(\sigma, 2\pi)$ monotonically increases in the range $0 \leq \sigma < \frac12$, and monotonically decreases in the range $\frac12 < \sigma \leq 1$.}

Setting $\sigma = \frac12 + \delta$, with $-\frac12 \leq \delta \leq \frac12$, when $t = 2\pi$, according to Eq. (\ref{gsigtdel}), we have:
    \beq\label{gsig2pi}
    G(\sigma, 2\pi) = 4\ln(\pi) - \sum^{\infty}_{n=0} \frac{\Psi\left(2n, \frac14 + \pi i \right) + \Psi\left(2n, \frac14 - \pi i \right)}{2^{2n-1}(2n)!}\delta^{2n}.
    \eeq
Neglecting terms with order higher than $O(\delta^6)$, we have:
    \beq\label{gsig2pidel}
    G(\sigma, 2\pi) \approx 4\ln(\pi) + g_0 + g_2\delta^2 + g_4\delta^4,
    \eeq
with $g_0 = -4.5746788$, $g_2 = -5.0986449 \times 10^{-2}$, $g_4 = 6.5574618 \times 10^{-4}$.

Using Eq. (\ref{gsig2pi}), we have:
    \beq
    \frac{dG(\sigma, 2\pi)}{d\sigma} = -\sum^{\infty}_{n=1}\frac{\Psi\left(2n, \frac14 + \pi i \right) + \Psi\left(2n, \frac14 - \pi i \right)}{2^{2n-1}(2n-1)!}\delta^{2n-1}.
    \eeq
{\it i.e.},
    \beq\label{dgdsigte}
    \frac{dG(\sigma, 2\pi)}{d\sigma} = g^{\prime}_1\delta + g^{\prime}_3\delta^3 + g^{\prime}_5\delta^5 + O(\delta^7),
    \eeq
where $g^{\prime}_{2k-1} = -\frac{\Psi\left(2n, \frac14 + \pi i\right) + \Psi\left(2n, \frac14 - \pi i\right)}{2^{2n-1}(2n -1)!}$, in particular:
$$g^{\prime}_1 = -1.0197290 \times 10^{-1},$$
$$g^{\prime}_3 = +2.6229847 \times 10^{-3},$$
and
$$g^{\prime}_5 = -6.8327078 \times 10^{-5}.$$
For a sufficiently large $N$, $|g^{\prime}_5| > \sum^N_{k=4}|g^{\prime}_{2k-1}| \rightarrow 1.904573728 \times 10^{-6}$, {\it i.e.}, neglecting the terms with order higher than $O(\delta^7)$ does not change the conclusion of $\frac{dG(\sigma, 2\pi)}{d\sigma}$ being positive or negative. 

According to Eq. (\ref{dgdsigte}), we have when $\delta < 0$: $g^{\prime}_1 \delta > 0$, $g^{\prime}_3 \delta^3 < 0$, $g^{\prime}_5\delta^5 > 0$, and $g^{\prime}_1 \delta + g^{\prime}_3 \delta^3 > 0$, so
$$
\frac{dG(\sigma, 2\pi)}{d\sigma} \approx g^{\prime}_1\delta + g^{\prime}_3\delta^3 + g^{\prime}_5\delta^5 > 0\hspace{0.25cm}\left(-\frac12 \leq \delta < 0 \right).
$$
When $\delta > 0$: $g^{\prime}_1 \delta < 0$, $g^{\prime}_3 \delta^3 > 0$, $g^{\prime}_5\delta^5 < 0$, and $g^{\prime}_1 \delta + g^{\prime}_3 \delta^3 < 0$, so
$$
\frac{dG(\sigma, 2\pi)}{d\sigma} \approx g^{\prime}_1\delta + g^{\prime}_3\delta^3 + g^{\prime}_5\delta^5 < 0\hspace{0.25cm}\left(0 < \delta \leq \frac12 \right).
$$
So $G(\sigma, 2\pi)$ monotonically increases in the range $0 \leq \sigma < \frac12$, and $G(\sigma, 2\pi)$ monotonically decreases in the range $\frac12 < \sigma \leq 1$.

{\bf 2) There are two symmetric zero points: $\sigma_{0\pm} = \frac12 + \delta_{\pm}$ of $G(\sigma, 2\pi)$ in the range $0 \leq \sigma < 1$}

According to Eq. (\ref{gsig2pidel}), setting $G(\sigma, 2\pi) = 0$, {\it i.e.},
$$
4\ln(\pi) + g_0 + g_2\delta^2 + g_4\delta^4 = 0,
$$
we have $\delta_{\pm} = \pm0.2885526325$.

{\bf 3) in the range $0 \leq \sigma < 1/2$, $|W(\sigma + 2\pi i)|$ monotonically decreases at the left side of $\sigma_{0-} = \frac12 + \delta_{-}$ and monotonically increases at the right side.}

According to the above item 1), $G(\sigma, 2 \pi)$ monotonically increases; according to the above item 2), $G\left( \frac12 + \delta_{-}, 2\pi \right) =0$, noticing that $G(0, 2\pi) = -0.008465084 < 0$,  $G\left(\frac12, 2\pi\right) = 0.004240720 > 0$, we find that during the continuous change of $\sigma$ from $0$ to $\frac12$, $G(\sigma, 2\pi)$ monotonically evolves from being less than zero to being larger than zero. Therefore, we have:
$$
\left.\frac{d|W(\sigma + 2\pi i)|}{d\sigma}\right|_{\sigma<\frac12+\delta_{-}} = \left. \frac{|W(\sigma + 2\pi i)|}4G(\sigma, 2\pi)\right|_{\sigma<\frac12+\delta_{-}} < 0,
$$
and
$$
\left.\frac{d|W(\sigma + 2\pi i)|}{d\sigma}\right|_{\sigma=\frac12+\delta_{-}} = \left. \frac{|W(\sigma + 2\pi i)|}4G(\sigma, 2\pi)\right|_{\sigma=\frac12+\delta_{-}} = 0,
$$
and
$$
\left.\frac{d|W(\sigma + 2\pi i)|}{d\sigma}\right|_{\sigma>\frac12+\delta_{-}} = \left. \frac{|W(\sigma + 2\pi i)|}4G(\sigma, 2\pi)\right|_{\sigma>\frac12+\delta_{-}} > 0.
$$
That is to say, $|W(\sigma + 2\pi i)|$ monotonically decreases at the left side of $\sigma_{0-} = \frac12 + \delta_{-}$ and monotonically increases at the right side.

{\bf 4) in the range $0 \leq \sigma < 1/2$, $|W(\sigma + 2\pi i)| < 1$.}

According to the above item 3), $|W(\sigma + 2\pi i)|$ monotonically decreases at the left side of $\sigma_{0-} = \frac12 + \delta_{-}$ and monotonically increases at the right side; therefore, the values of $|W(\sigma + 2\pi i)|$ at the two ends: $\sigma = 0$ and $\sigma = \frac12$ must be larger than its value at points between the two ends. Because $|W(0 + 2\pi i)| = 0.9999999991 < 1$, and $|W(1/2 + 2\pi i)| = 1$, we have $|W(\sigma_{-} + 2\pi i)| < 1$ for $\left( 0 \leq \sigma < \frac12 \right)$.

\section{When $|\zeta(s)| = |\zeta(1 - s)|$, it must be true that $|W(s)| = 1$}\label{equiv}

{\it Proof}

According to the Riemann-Zeta equation $\zeta(s) = W(s) \zeta(1 - s)$, we have:
$    W(s) = \frac{\zeta(s)}{\zeta(1-s)}$.
Because ${\overline{W}}(s) = W(s^{\ast})$, we have $\frac{{\overline{\zeta}}(s)}{{\overline{\zeta}}(1-s)} = \frac{\zeta(s^{\ast})}{\zeta(1-s^{\ast})}$.
Therefore:
    \beq\label{wsq}
    |W(s)|^2 = W(s) {\overline{W}}(s) = \frac{\zeta(s){\overline{\zeta}}(s)}{\zeta(1-s){\overline{\zeta}}(1-s)} = \frac{\zeta(s) \zeta\left(s^{\ast}\right)}{\zeta(1 - s) \zeta\left( 1 - s^{\ast} \right)}.
    \eeq
Here, $s = \sigma + it$, and $s^{\ast} = \sigma - it$.

Assuming
    \beq\label{zetaeq}
    |\zeta(s)| = |\zeta(1-s)|,
    \eeq
we have:
    \beq\label{zetazeta}
    |\zeta(s)|^2 = |\zeta(1-s)|^2,
    \eeq
where $|\zeta(s)|^2 = \zeta(s) \zeta(s^{\ast})$ and $|\zeta(1-s)|^2 = \zeta(1-s) \zeta(1-s^{\ast})$ are both real function of $\sigma$ and $t$. For the convenience of discussion, let us denote:
$$P(\sigma, t) \equiv \zeta(s)\zeta(s^{\ast}),$$
and
$$Q(\sigma, t) \equiv \zeta(1-s)\zeta(1-s^{\ast}).$$

1) If $\zeta(s) \neq 0$, then it must be true that $\zeta(1-s) \neq 0$. According to Eq. (\ref{wsq}) and Eq. (\ref{zetaeq}), it must be true that:
    \beq
    |W(s)| = \frac{|\zeta(s)|}{|\zeta(1-s)|} = 1.
    \eeq

2) If $\zeta(s) = 0$, then it must be true that $\zeta(1-s) = 0$. In this case, we have $|W(s)| = \frac00$. In this following, we prove the limit of this $\frac00$ is $1$.

2.1) In the range $0 \leq \sigma \leq 1$, and $t \neq 0$, due to the fact that $\zeta(s)$, $\zeta(s^{\ast})$, $\zeta(1-s)$, and $\zeta(1-s^{\ast})$ are all analytical functions. Furthermore, there exists their arbitrary order of partial derivative with respect to $\sigma$ and $t$: $\frac{\partial^m\zeta(s)}{\partial \sigma^m}$, $\frac{\partial^m\zeta(s)}{\partial t^m}$, $\frac{\partial^n\zeta(s^{\ast})}{\partial \sigma^n}$, $\frac{\partial^n\zeta(s^{\ast})}{\partial t^n}$; $\frac{\partial^m\zeta(1-s)}{\partial \sigma^m}$, $\frac{\partial^m\zeta(1-s)}{\partial t^m}$, $\frac{\partial^n\zeta(1-s^{\ast})}{\partial \sigma^n}$, $\frac{\partial^n\zeta(1-s^{\ast})}{\partial t^n}$. Therefore the real functions $P(\sigma, t)$ and $Q(\sigma, t)$ are partial differentiable to arbitrary order with respect to $\sigma$ and $t$.

2.2) If $\zeta(s)$ has zeros: $s_n = \sigma_n + it_n$ with ($n=1, 2, 3, \cdots$), so that $\zeta(s_n) = 0$. Then the complex conjugate ${\overline{\zeta}}(s_n) = \zeta(s^{\ast}_n) = 0$. Therefore, for these $\sigma_n$ and $t_n$, they are the zeros of the real function $P(\sigma_n, t)$. Then there exists the Taylor expansion of $P(\sigma_n, t)$ around the zeros, $t_n$:
    \beq
    P(\sigma_n, t) =\left. \sum^{\infty}_{j = 0} \frac1{j!}\frac{d^j P(\sigma_n, t)}{dt^j}\right|_{t=t_n}(t-t_n)^j
    \eeq
where the derivatives in the expansion coefficients: $\left\{ \left. \frac{d^j P(\sigma_n, t)}{dt^j}\right|_{t = t_n}, j =0,1,2,\cdots,\infty \right\}$ won't be all equal to zero according to the series expansion theorem. Notice that, the zeroth order expansion coefficient of $P(\sigma_n, t)$ is just $P(\sigma_n, t) \equiv \zeta(s_n) \zeta(s^{\ast}_n) = 0$. The first order expansion coefficient is then:
$$
\left. \frac{dP(\sigma_n, t)}{dt}\right|_{t = t_n} = \left. \left( \frac{d\zeta(\sigma_n + it)}{dt}\zeta(\sigma_n - it) + \zeta(\sigma_n + it) \frac{d\zeta(\sigma_n - it)}{dt} \right)\right|_{t=t_n} = 0.
$$

There must exit a certain $k \geq 2$, so that:
    \beq\label{dPdtk}
    \lim_{t \rightarrow t_n} \left( \frac{d^k P(\sigma_n, t)}{dt^k}\right) = \lim_{t\rightarrow t_n} \left( \left. \frac{d^k}{dt^k}\sum^{\infty}_{j = 0} \frac1{j!}\frac{d^j P(\sigma_n, t)}{dt^j}\right|_{t = t_n}(t - t_n)^j\right) = \left. \frac{d^k P(\sigma_n, t)}{dt^k} \right|_{t = t_n} \neq 0.
    \eeq

2.3) With the assumption in Eq. (\ref{zetaeq}), following Eq. (\ref{zetazeta}), it must be true that: $|\zeta(1 - s_n)|^2 = 0$, {\it i.e.}, $t_n$ is the zeros of the real function $Q(\sigma_n, t)$. Similarly, there exists the Taylor expansion around the zeros $t_n$:
    \beq
    Q(\sigma_n, t) = \left. \sum^{\infty}_{j = 0} \frac1{j!} \frac{d^j Q(\sigma_n, t)}{dt^j} \right|_{t = t_n}(t - t_n)^j.
    \eeq
where the derivatives in the expansion coefficients: $\left\{ \left. \frac{d^j Q(\sigma_n, t)}{dt^j}\right|_{t = t_n}, j =0,1,2,\cdots,\infty \right\}$ won't be all equal to zero. Because both the zeroth order expansion coefficient and the first order expansion coefficient are equal to zero; there must exit a certain order $k \geq 2$, so that:
    \beq\label{dQdtk}
    \lim_{t \rightarrow t_n} \left( \frac{d^k Q(\sigma_n, t)}{dt^k}\right) = \lim_{t\rightarrow t_n} \left( \left. \frac{d^k}{dt^k}\sum^{\infty}_{j = 0} \frac 1{j!}\frac{d^j Q(\sigma_n, t)} {dt^j} \right|_{t = t_n} (t - t_n)^j \right) = \left. \frac{d^k Q(\sigma_n, t)}{dt^k} \right|_{t = t_n} \neq 0.
    \eeq
According to Eq. (\ref{zetazeta}), Eq. (\ref{dPdtk}), and Eq. (\ref{dQdtk}), there must exist a certain $k \geq 2$, so that:
    \beq\label{dkneq0}
    \left. \frac{d^k P(\sigma_n, t)}{dt^k}\right|_{t = t_n} = \left. \frac{d^k Q(\sigma_n, t)}{dt^k}\right|_{t = t_n} \neq 0.
    \eeq
For $|W(\sigma_n + it_n)|^2 = \frac{P(\sigma_n, t)}{Q(\sigma_n, t)} = \frac00$, because:
    \beq\label{limW2}
    \lim_{t \rightarrow t_n} |W(\sigma_n + it)|^2 = \lim_{t \rightarrow t_n} \frac{P(\sigma_n, t)}{Q(\sigma_n, t)} = \lim_{t \rightarrow t_n} \frac{\frac{d^k}{dt^k}P(\sigma_n, t)}{\frac{d^k}{dt^k}Q(\sigma_n, t)} = \frac{\left. \frac{d^kP(\sigma_n, t)}{dt^k}\right|_{t=t_n}}{\left.\frac{d^kQ(\sigma_n, t)}{dt^k}\right|_{t=t_n}}.
    \eeq
According to Eq. (\ref{dkneq0}) and Eq. (\ref{limW2}), we have:
    \beq\label{limW2eq1}
    \lim_{t\rightarrow t_n} |W(\sigma_n + it)|^2 = 1.
    \eeq
Based on the Taylor expansion of the real functions $P(\sigma_n, t)$ and $Q(\sigma_n, t)$ around their zeros, we can prove by the same procedure that there exists a certain $k \geq 2$, so that:
    \beq
    \left. \frac{d^k P(\sigma, t_n)}{d\sigma^k}\right|_{\sigma = \sigma_n} = \left. \frac{d^k Q(\sigma, t_n)}{d\sigma^k}\right|_{\sigma = \sigma_n} \neq 0.
    \eeq
Therefore, we have:
    \beq\label{limW2sig}
    \lim_{\sigma \rightarrow \sigma_n} |W(\sigma + it_n)|^2 = \lim_{\sigma \rightarrow \sigma_n} \frac{P(\sigma, t_n)}{Q(\sigma, t_n)} = \lim_{\sigma \rightarrow \sigma_n} \frac{\frac{d^k}{d\sigma^k}P(\sigma, t_n)}{\frac{d^k}{d\sigma^k}Q(\sigma, t_n)} = \frac{\left. \frac{d^k P(\sigma, t_n)}{d\sigma^k}\right|_{\sigma=\sigma_n}}{\left.\frac{d^kQ(\sigma, t_n)}{d\sigma^k}\right|_{\sigma=\sigma_n}} = 1.
    \eeq
According to Eq. (\ref{limW2eq1}) and Eq. (\ref{limW2sig}), we have:
    \beq
    \lim_{t \rightarrow t_n} |W(\sigma_n + it)| = 1,
    \eeq
and
    \beq
    \lim_{\sigma \rightarrow \sigma_n} |W(\sigma + it_n)| = 1.
    \eeq
The above serves as a proof.

\section{Proof of $\frac d{dt}\left| \Gamma \left( \frac{1-s}2\right)\right| = \sum^{\infty}_{n=1}\frac{-t\left| \Gamma\left( \frac{1-s}2\right)\right|}{|it + 2n -\sigma - 1|^2}$, $\frac d{dt}\left| \Gamma \left( \frac{s}2\right)\right| = \sum^{\infty}_{n=1}\frac{-t\left| \Gamma\left( \frac{s}2\right)\right|}{|it + 2n + \sigma - 2|^2}$}\label{dgammadt}

Because:
    \beq\label{gammas2}
    \left| \Gamma\left( \frac s 2\right)\right|^2 = \Gamma\left(\frac s2\right) \Gamma\left(\frac{s^{\ast}}2\right).
    \eeq
    \beq\label{gamma1s2}
    \left| \Gamma\left( \frac {1-s} 2\right)\right|^2 = \Gamma\left(\frac {1-s}2\right) \Gamma\left(\frac{1- s^{\ast}}2\right).
    \eeq
Furthermore, in the range $0 \leq \sigma \leq 1$ and $t \neq 0$, both $\Gamma\left(\frac s2 \right)$ and $\Gamma\left(\frac{1-s}2\right)$ are analytical functions.

We take derivative of Eq. (\ref{gammas2}) and Eq. (\ref{gamma1s2}) with respect to $t$:
    \beq\label{dgamma2dt}
    \frac d{dt}\left| \Gamma\left( \frac s 2\right)\right|^2 = \frac i2 \left| \Gamma\left( \frac s 2\right)\right|^2 \left[\Psi\left(\frac s2\right) - \Psi\left(\frac{s^{\ast}}2\right)\right],
    \eeq
    \beq\label{dgamma1ms2dt}
    \frac d{dt}\left| \Gamma\left( \frac {1-s} 2\right)\right|^2 = \frac i2 \left| \Gamma\left( \frac {1-s} 2\right)\right|^2 \left[\Psi\left(\frac {1-s^{\ast}}2\right) - \Psi\left(\frac{1 - s}2\right)\right].
    \eeq

According to Eq. (\ref{Psigamma}),
    \beq\label{PsimPsi}
    \Psi\left(\frac{\sigma - it}2\right) - \Psi\left(\frac{\sigma + it}2\right) = \sum^{\infty}_{n=1}\frac{-4it}{|it+ 2n + \sigma -2|^2},
    \eeq
    \beq\label{nPsipPsi}
    -\Psi\left(\frac{1-\sigma - it}2\right) + \Psi\left(\frac{1 - \sigma + it}2\right) = \sum^{\infty}_{n=1}\frac{4it}{|it+ 2n - \sigma -1|^2},
    \eeq
Plugging Eq. (\ref{PsimPsi}) into Eq. (\ref{dgamma2dt}), and plugging Eq. (\ref{nPsipPsi}) into Eq. (\ref{dgamma1ms2dt}), we have:
    \beq\label{dgamma2dtseries}
    \frac d{dt}\left| \Gamma\left( \frac s 2\right)\right|^2 = \frac i2 \left| \Gamma\left( \frac s 2\right)\right|^2 \sum^{\infty}_{n=1}\frac{4it}{|it+ 2n + \sigma -2|^2}.
    \eeq
    \beq\label{dgamma1ms2dtseries}
    \frac d{dt}\left| \Gamma\left( \frac {1-s} 2\right)\right|^2 = \frac i2 \left| \Gamma\left( \frac {1-s} 2\right)\right|^2 \sum^{\infty}_{n=1}\frac{4it}{|it+ 2n - \sigma -1|^2}.
    \eeq

Notice that: $\frac d{dt}\left| \Gamma\left( \frac s 2\right)\right|^2 = 2 \left| \Gamma\left( \frac s 2\right)\right| \frac d{dt}\left| \Gamma\left( \frac s 2\right)\right|$, $\frac d{dt}\left| \Gamma\left( \frac {1-s} 2\right)\right|^2 = 2 \left| \Gamma\left( \frac {1-s} 2\right)\right| \frac d{dt}\left| \Gamma\left( \frac {1-s} 2\right)\right|$, therefore, Eq. (\ref{dgamma2dtseries}) and Eq. (\ref{dgamma1ms2dtseries}) can finally simplified as:
    \beq
    \frac d{dt}\left| \Gamma\left( \frac s 2\right)\right| = -t \sum^{\infty}_{n=1}\frac{\left| \Gamma\left( \frac s 2\right)\right|}{|it+ 2n + \sigma -2|^2},
    \nonumber
    \eeq
    \beq
    \frac d{dt}\left| \Gamma\left( \frac {1 - s} 2\right)\right| = -t \sum^{\infty}_{n=1}\frac{\left| \Gamma\left( \frac {1 - s} 2\right)\right|}{|it+ 2n - \sigma -1|^2}.
    \nonumber
    \eeq

This serves as a proof.

\section{For an arbitrary finite $t > 0$, the $\frac{|\zeta(s)|}{|\zeta(1-s)|}$ satisfying $\frac d{dt} \left( \frac{|\zeta(s)|}{|\zeta(1-s)|}\right) = \frac d{dt} |W(s)| \neq 0$, can not be in the form of $\frac00$.}\label{appendzetazetano00}

{\it Proof}

According to Eq. (\ref{dwdtpos}), when $|W(s)| >0$, and $\sigma < \frac12$, we always have:
    \beq
    \frac d{dt} \left( \frac{|\zeta(s)|}{|\zeta(1-s)|}\right) = \frac{|\zeta(1-s)|\frac{d|\zeta(s)|}{dt} - |\zeta(s)|\frac{d|\zeta(1-s)|}{dt}}{|\zeta(1-s)|^2} = \frac d{dt} |W(s)| > 0
    \eeq
therefore:
    \beq\label{zdzg0}
    |\zeta(1-s)|\frac{d|\zeta(s)|}{dt} - |\zeta(s)|\frac{d|\zeta(1-s)|}{dt} > 0.
    \eeq

Notice that:
$$
\frac{d|\zeta(s)|^2}{dt} = 2|\zeta(s)|\frac{d|\zeta(s)|}{dt} = \frac{d\zeta(s)}{dt}\zeta\left(s^{\ast}\right) +  \frac{d\zeta\left(s^{\ast}\right)}{dt}\zeta(s),
$$
we have:
$$
\frac{d|\zeta(s)|}{dt} = \frac1{2|\zeta(s)|}\left[ \frac{d\zeta(s)}{dt}\zeta\left(s^{\ast}\right) +  \frac{d\zeta\left(s^{\ast}\right)}{dt}\zeta(s) \right].
$$

Similarly, we have:
$$
\frac{d|\zeta(1- s)|}{dt} = \frac1{2|\zeta(1-s)|}\left[ \frac{d\zeta(1-s)}{dt}\zeta\left(1-s^{\ast}\right) +  \frac{d\zeta\left(1-s^{\ast}\right)}{dt}\zeta(1-s) \right].
$$

Therefore, Eq. (\ref{zdzg0}) can be rewritten as:
    \beq\label{halfg0}
    \frac12\left\{\frac{\frac{d\zeta(s)}{dt}\zeta\left(s^{\ast}\right) +  \frac{d\zeta\left(s^{\ast}\right)}{dt}\zeta(s)}{|W(s)|} - |W(s)| \left[ \frac{d\zeta(1-s)}{dt}\zeta\left(1-s^{\ast}\right) +  \frac{d\zeta\left(1-s^{\ast}\right)}{dt}\zeta(1-s) \right]\right\} > 0.
    \eeq

Notice that, for an arbitrary finite $t > 0$, $|W(s)| > 0$ and have no singular points, $\zeta(s)$, $\zeta\left(s^{\ast}\right)$, $\zeta(1-s)$, $\zeta\left(1-s^{\ast}\right)$ are all analytical functions. Furthermore, their derivatives $\frac{d\zeta(s)}{dt}$, $\frac{d\zeta\left(s^{\ast}\right)}{dt}$, $\frac{d\zeta(1-s)}{dt}$, $\frac{d\zeta\left(1-s^{\ast}\right)}{dt}$ do not have singular points. If there exists a certain $s_n$, so that $|\zeta\left(s_n\right)| = 0$, $|\zeta\left(1 - s_n\right)| = 0$, then it must be true that $\zeta\left(s_n\right) = \zeta\left(s_n^{\ast}\right) = 0$, $\zeta\left(1 - s_n\right) = \zeta\left( 1 - s_n^{\ast}\right) = 0$. Plugging these into Eq. (\ref{halfg0}), we have:
$$
    \frac12\left\{ \frac{0\cdot \left.\frac{d\zeta(s)}{dt}\right|_{s=s_n} +  0 \cdot \left. \frac{d\zeta\left(s^{\ast}\right)}{dt}\right|_{s^{\ast} = s_n^{\ast}}}{|W(s_n)|} - |W(s_n)| \left[ 0\cdot\left.\frac{d\zeta(1-s)}{dt}\right|_{s = s_n} +  0\cdot \left. \frac{d\zeta\left(1-s^{\ast}\right)}{dt}\right|_{s^{\ast} = s_n^{\ast}} \right]\right\} = 0 > 0.
$$

This is to say that $\left|\zeta\left(s_n\right)\right| = 0$, $\left|\zeta\left(1-s_n\right)\right| = 0$ must lead to the invalidation of inequality relation in Eq. (\ref{halfg0}). Therefore, for $\sigma < \frac 12$, and an arbitrary finite $t > 0$, the $\frac{|\zeta(s)|}{|\zeta(1-s)|}$ satisfying $\frac d{dt} \left( \frac{|\zeta(s)|}{|\zeta(1-s)|}\right) > 0$, can not be in the form of $\frac00$.

Similarly, for $\sigma > \frac 12$, and an arbitrary finite $t > 0$, the $\frac{|\zeta(s)|}{|\zeta(1-s)|}$ satisfying $\frac d{dt} \left( \frac{|\zeta(s)|}{|\zeta(1-s)|}\right) < 0$, can not be in the form of $\frac00$.

This serves as a proof.

\nocite{*}

\section*{References}

\end{document}